\documentclass[11pt, a4paper]{article}
\usepackage[top=2.7cm,bottom=2.7cm,left=2.5cm,right=2.5cm]{geometry}

\usepackage{graphicx,caption,amsmath,amsfonts,amssymb,booktabs}
\usepackage{float,listings}
\floatstyle{ruled}
\newfloat{algorithm}{htbp}{loa}
\floatname{algorithm}{Algorithm}
\usepackage[dvipsnames]{xcolor}
\PassOptionsToPackage{hyphens}{url}
\usepackage[numbers,sort&compress]{natbib}
\usepackage[breaklinks,colorlinks,allcolors=black]{hyperref}
\DeclareTextSymbolDefault{\DJ}{T1}
\DeclareTextSymbolDefault{\dj}{T1}

\definecolor{darkgreen}{rgb}{0.0, 0.5, 0.0} 
\definecolor{darkred}{rgb}{0.5, 0.0, 0.0}

\lstdefinestyle{fortran-style}{
    language=[90]Fortran,
    basicstyle=\ttfamily\small,
    morekeywords={abstract, import, asinh, acosh , atanh},
    keywordstyle=\color{blue}\bfseries,
    commentstyle=\color{darkgreen},
    stringstyle=\color{darkred},
    numbers=left,
    numberstyle=\tiny\color{gray},
    stepnumber=1,
    numbersep=5pt,
    showspaces=false,
    showstringspaces=false,
    showtabs=false,
    tabsize=4,
    captionpos=b,
    breaklines=true,
    breakatwhitespace=false,
    escapeinside={\%*}{*)}
}

\lstdefinestyle{python-style}{
    language=Python,
    basicstyle=\ttfamily\small,
    keywordstyle=\color{blue},
    commentstyle=\color{darkgreen},
    stringstyle=\color{darkred},
    showstringspaces=false,
    keywords=[2]{grad_k[k].item()},
    keywordstyle=[2]{\color{black}},
    numbers=left,
    numberstyle=\tiny\color{gray},
    breaklines=true,
    captionpos=b
}

\lstdefinelanguage{Julia}{
    keywords={struct, abstract, type, mutable, function, end, where, if, 
    else, elseif, for, while, in, return, break, continue, try, catch, 
    finally, local, global, const, let, do, begin, quote, module, 
    baremodule, include},
    sensitive=true,
    morecomment=[l]{\#},
    morestring=[b]",
}

\lstdefinestyle{julia-style}{
    language=Julia,
    basicstyle=\ttfamily\small,
    morekeywords={import, using, function, end, mutable, struct, 
    abstract, type, where, if, else, elseif, for, while, in, return, 
    break, continue, try, catch, finally, local, global, const, let, do, 
    begin, quote, module, baremodule, include},
    keywordstyle=\color{blue}\bfseries,
    commentstyle=\color{darkgreen},
    stringstyle=\color{darkred},
    numbers=left,
    numberstyle=\tiny\color{gray},
    stepnumber=1,
    numbersep=5pt,
    showspaces=false,
    showstringspaces=false,
    showtabs=false,
    tabsize=4,
    captionpos=b,
    breaklines=true,
    breakatwhitespace=false,
    escapeinside={\%*}{*)}
}

\title{Dual Numbers for Arbitrary Order Automatic Differentiation}

\author{
  F. Pe\~nu\~nuri,~
  K. B. Cant\'un-Avila,~
  R. Pe\'on-Escalante
}

\date{}

\begin{document}

\maketitle

\vspace*{-0.7cm}
\begin{center}
\footnotesize{Facultad de Ingenier\'ia, Universidad
Aut\'onoma de Yucat\'an, A.P. 150, Cordemex, M\'erida, Yucat\'an,
M\'exico.\\
Corresponding author: F. Pe\~nu\~nuri
(\href{mailto:francisco.pa@correo.uady.mx}{\texttt{francisco.pa@correo.uady.mx}}).}
\end{center}

\vspace*{0.5cm}
\begin{abstract}
We present DNAOAD-FPM, a modern Fortran package for univariate
automatic differentiation with arbitrary-order dual numbers. Its
complex-valued components store ordinary derivatives through a
user-selected order within a compiled capacity. A generic composition
routine computes all requested orders using one table of partial Bell
polynomials per composition. For $\exp(\sin x)$ at $x=0.5$,
DNAOAD-FPM had lower median evaluation times at the tested orders
through 8, and ADOL-C's higher-order forward sweep with a reused tape
had lower medians at the tested orders from 12 to 24. Two applications illustrate how dual results pass between routines. In the first, a Newton--Raphson solver for Kepler's
equation returns the dual eccentric anomaly to a routine that computes
orbital position. Checking every active residual component prevents
premature termination at a converged scalar root. An eighth-order
expansion met a position tolerance of $10^{-12}$ times the semimajor
axis at all 256 query points at each of two centers. Including
construction, the ratio of direct scalar solution time to expansion
time ranged from 1.03 to 3.33 across the two centers and two
compilers. In the second application, a geometry routine returns an
orthonormal frame attached to the coupler of a spherical four-bar
linkage. A single evaluation of this frame through order eight
supplies the coupler's angular velocity and its first seven time
derivatives. All components through angular snap reproduce published
values to their reported precision. At the tested configuration, the
four higher-order vectors agree with an independent calculation to
within $1.2\times10^{-15}$ in relative Euclidean norm.
\end{abstract}


\vspace*{0.7cm}
\section{Introduction}
Derivatives enter optimization algorithms, sensitivity analyses, and
numerical models throughout science and engineering. Symbolic
differentiation can produce expressions that grow rapidly with
derivative order, while finite-difference approximations require a
choice of step size that balances truncation and roundoff errors.
Automatic differentiation (AD) computes derivatives of functions
defined by numerical programs by applying differentiation rules to the
operations performed during evaluation
\cite{Griewank1989,Bucker2006,Atilim2018}. Many AD tools address
first- and second-order derivatives of real-valued functions
\cite{Phipps2012,Quang2012,Hascoet2013,Hogan2014,Weinstein2016,
  Felix2016,Slusanschi2016}. We consider programs in which the result
of one calculation becomes the input to another. In DNAOAD-FPM, a
Fortran routine can return a value together with its derivatives
through a chosen order in one independent variable. Subsequent
routines can use that object directly. A root solver, for example, can
return a dual root to a routine that computes a physical observable.

One early application of automatic Taylor coefficient generation
was the numerical integration of ordinary differential equations.
In 1971, Barton, Willers, and Zahar described how to generate a
program that advances the solution of a system of ordinary
differential equations by Taylor series \cite{Barton1971}. Chang and
Corliss later described how ATOMFT recursively generated the series
coefficients by automatic differentiation \cite{Chang1994}. By 1996,
ADOL-C provided higher-order differentiation for functions defined by
C or C++ programs \cite{Griewank1996}. Its higher-order forward mode
propagates univariate Taylor coefficients through a recorded sequence
of operations \cite{Walther2009,adolcgithub}.

In univariate Taylor mode, an intermediate quantity $u(t)$ is
represented through order $N$ by its normalized coefficients
$u_k=u^{(k)}(t_0)/k!$. For many elementary functions, a differential
identity gives a recurrence for the output coefficients. For example,
if $y=\exp(u)$, then $y'=u'y$ determines each successive coefficient
of $y$ from the coefficients of $u$ and the preceding coefficients
of $y$. Through order $N$, this recurrence requires $O(N^2)$
arithmetic operations and $O(N)$ storage \cite{Griewank2008}.
Truncated Taylor arithmetic is also available through overloaded
operators in the Fortran~95 module TaylUR
\cite{vonHippel2006,vonHippel2010}, in AuDi and its Python interface
pyaudi under the name generalized dual numbers
\cite{Izzo2017,Izzo2026}, and in TaylorSeries.jl for Julia
\cite{Benet2019}. Taylor-series integration is also implemented in
TIDES \cite{Abad2012}.

The higher-order dual numbers used in DNAOAD-FPM belong to this same
truncated polynomial algebra. With a single nilpotent generator
satisfying $\epsilon^{N+1}=0$, the choice of basis
$\epsilon_k=\epsilon^k/k!$ gives
\[
  \sum_{k=0}^{N}\frac{u^{(k)}(t_0)}{k!}\,\epsilon^k
  = \sum_{k=0}^{N}u^{(k)}(t_0)\,\epsilon_k .
\]
Thus the $N+1$ components in the divided-power basis are the ordinary
derivatives themselves. Messelmi studied the algebraic properties of
this extension \cite{Messelmi2015}. Szirmay-Kalos used the same basis
for generalized dual numbers and gave \texttt{C++} implementations
of elementary functions through second order \cite{Kalos2021}.
In previous work, we used dual numbers through fourth order to
compute kinematic quantities and directional derivatives
\cite{Peon2024,Penunuri2024}.

Higher derivatives can also be obtained by repeated differentiation.
In graph-based frameworks such as PyTorch, this can involve
differentiating graphs produced by earlier derivative evaluations
\cite{paszke2019pytorch}. ForwardDiff in Julia supports recursive
nesting of first-order dual numbers \cite{Revels2016}. Nesting $N$
scalar first-order dual levels gives a representation with $2^N$
scalar components. The multidual construction in MultiZ likewise
uses a tensor-product algebra with distinct nilpotent generators
\cite{AguirreMesa2020}. These representations retain information
about combinations of differentiation directions. For derivatives
along a single direction, the $N+1$ components above suffice.

Within the univariate algebra, a computational choice remains in
how to evaluate the composition $f(u)$ through every order up to $N$.
TaylUR uses Faà di Bruno's formula and, in its original
implementation, generates integer partitions and stores their weights
\cite{vonHippel2006}. The implementation studied here evaluates
partial Bell polynomials by dynamic programming. For each elementary
function handled by the generic composition routine, one Bell table
serves all requested derivative orders. We analyze the work and
storage required by this calculation. For $h(x)=\exp(\sin x)$, we
compare evaluation times and numerical errors with a Taylor recurrence
specialized to $h$ and with the higher-order forward mode of ADOL-C.

DNAOAD-FPM is distributed through its Fortran Package Manager (FPM)
repository \cite{Penunuri_DNAOAD_FPM}. We use two applications to
examine how the dual type works across routine boundaries. In the
first, a Newton--Raphson solver returns the eccentric anomaly in
Kepler's equation. Its derivatives determine local orbital positions,
so we test both the solver's stopping criterion and the cost of
predicting several positions from a single expansion. In the second,
a geometry routine returns a frame attached to the coupler of a
spherical four-bar linkage. From that frame, the angular velocity
expression gives the coupler's angular velocity and its first seven
time derivatives.

\section{Dual numbers and derivatives}

\subsection{First-order case}

A first-order dual number has the form $g=g_0+g_1\epsilon_1$,
where the coefficients may be real or complex and $\epsilon_1^2=0$.
For a function $f$ analytic at $g_0$, its Taylor expansion in the
nilpotent increment gives
\begin{equation}
  f(g_0+g_1\epsilon_1)=f(g_0)+f'(g_0)g_1\epsilon_1 .
\end{equation}
For example, $\sin(g)=\sin(g_0)+g_1\cos(g_0)\epsilon_1$.
Taking $g_1=1$ returns the ordinary derivative of the sine at $g_0$.
When $g$ comes from an earlier calculation, its component $g_1$
carries the derivative of that intermediate result, and the same
formula applies the chain rule.

First- and second-order constructions have been used in a range of
applications \cite{Cheng1994,Jefrey2011,Wenbin2013,Penunuri2020,Kalos2021}.
To propagate higher orders in a single direction, DNAOAD-FPM stores
one component per derivative order in a flat array.

\subsection{Arbitrary-order case}
\label{sec:arbitrary}

Derivatives of arbitrary order can be computed by defining a dual number
of order $n$ as \cite{Kalos2021,Penunuri2024}
\begin{align}
\label{dualn}
r_n = \sum_{k=0}^{n} a_k\,\epsilon_k,
\end{align}
where $a_k$ are complex coefficients and the basis elements
$\epsilon_k$ satisfy the multiplication rule
\begin{equation}
\label{gentabmult}
\epsilon_i \cdot \epsilon_j =
\begin{cases}
0, & \text{if } i+j > n, \\[0.2cm]
\dfrac{(i+j)!}{i!\,j!}\,\epsilon_{i+j}, & \text{otherwise}.
\end{cases}
\end{equation}
Writing $\epsilon_k = \epsilon^k/k!$ with $\epsilon^{n+1}=0$,
Eq.~(\ref{gentabmult}) follows from $\epsilon^i\epsilon^j=\epsilon^{i+j}$.
The algebra is therefore that of the multidual numbers of
\cite{Messelmi2015}, that is, of truncated Taylor polynomials, written in
a basis in which the components are the derivatives themselves
\cite{Kalos2021}.

Evaluating the Taylor expansion of an analytic function at
$z+\epsilon_1$ and using Eq.~(\ref{gentabmult}) yields
\begin{align}
\label{fdnfz}
f(z+\epsilon_1)
= f(z)\,\epsilon_0 + f'(z)\,\epsilon_1 + \cdots + f^{(n)}(z)\,\epsilon_n.
\end{align}
For a general inner argument $g$, each component of $f(g)$ depends
on derivatives of both functions. The following section gives the
composition rule used by the package \cite{Johnson2002,Craik2005}.

\section{Implementing arbitrary-order dual arithmetic}
Listing~\ref{lst:dualzn} shows the Fortran type used for scalar values,
intermediate results, and outputs.

\begin{lstlisting}[
style=fortran-style,
caption={Fortran derived type \texttt{dualzn} using the precision
specified by \texttt{prec}. The $k$-th coefficient of a dual number $g$
is accessed as \texttt{g\%f(k)}.},
label={lst:dualzn}]
type, public :: dualzn
   complex(prec), dimension(0:max_order_dualzn) :: f
   integer :: ord
end type dualzn
\end{lstlisting}

Each object stores $M+1$ complex values and an integer active order
$N=\texttt{ord}$, where $M=\texttt{max\_order\_dualzn}$ is selected at
compilation. Only components $0$ through $N$ are active, with
$0\le N\le M$. Increasing the capacity requires recompilation.

\subsection{Chain rule via Fa\`a di Bruno formula}

According to the Fa\`a di Bruno formula, the $n$-th derivative of a
composite function $f(g(x))$ is given by
\begin{align}
\label{FaDiB}
D^n f(g(x)) =
\sum \frac{n!}{k_1! k_2! \cdots k_n!}
f^{(k_1 + k_2 + \cdots + k_n)}(g(x))
\prod_{j=1}^n
\left( \frac{g^{(j)}(x)}{j!} \right)^{k_j},
\end{align}
where $D^n = d^n/dx^n$, $f^{(k)}(x) = D^k f(x)$, and the sum is taken over
all nonnegative integer solutions of the Diophantine equation
\begin{align}
k_1 + 2k_2 + 3k_3 + \cdots + n k_n = n.
\end{align}

Although recursive formulations that avoid explicitly solving this
Diophantine equation exist \cite{Mohammed2016}, an equivalent formulation is obtained by rewriting
Eq.~(\ref{FaDiB}) as
\begin{align}
\label{DnChainRB}
D^n f(g(x)) =
\sum_{k=1}^n
f^{(k)}(g(x))\,
B_{n,k}\!\left(g'(x), g''(x), \dots, g^{(n-k+1)}(x)\right),
\end{align}
where $B_{n,k}$ are the partial Bell polynomials.

The partial Bell polynomials can be defined recursively as
\cite{Birmajer2015}
\begin{align}
B_{n,k}(x_1,\dots,x_{n-k+1})
&= \sum_{i=0}^{n-k}
\binom{n-1}{i}\,
x_{i+1}\,
B_{n-i-1,k-1}(x_1,\dots,x_{n-k-i+1}),
\end{align}
with initial conditions
\begin{align}
B_{0,0} &= 1, \\
B_{n,0} &= 0 \quad (n \ge 1), \\
B_{0,k} &= 0 \quad (k \ge 1).
\end{align}

\subsection{Computing a complete composition}
The private routine \texttt{compose\_dual} requests the outer derivative
vector once. It evaluates Eq.~(\ref{DnChainRB}) for all orders using one
shared Bell table. Algorithm~\ref{BellY} uses $g_j=g^{(j)}(x_0)$ and
$F_k=f^{(k)}(g_0)$. The Fortran arrays for the binomial coefficients and
Bell polynomials are rectangular, with triangular index ranges.

\begin{algorithm}
\caption{Composition through the active order $N$.}\label{BellY}
\begin{enumerate}\setlength{\itemsep}{0pt}\setlength{\parskip}{0pt}
\item Compute $F_0,\ldots,F_N$ at $g_0$ and set $y_0=F_0$
\item \textbf{if} $N=0$ \textbf{then return} $y$
\item Fill $C_{r,j}=\binom{r}{j}$ by Pascal's recurrence for $0\le j\le r<N$
\item Set $B_{0,0}=1$ and $B_{n,0}=0$ for $1\le n\le N$
\item \textbf{for} $n=1$ to $N$ \textbf{do}
\item \hspace*{1em}$y_n\gets0$
\item \hspace*{1em}\textbf{for} $k=1$ to $n$ \textbf{do}
\item \hspace*{2em}$B_{n,k}\gets\sum_{j=1}^{n-k+1} C_{n-1,j-1}g_j B_{n-j,k-1}$
\item \hspace*{2em}$y_n\gets y_n+F_k B_{n,k}$
\item \hspace*{1em}\textbf{end for}
\item \textbf{end for}
\item \textbf{return} $(y_0,\ldots,y_N)$
\end{enumerate}
\end{algorithm}

The inner sums comprise $N(N+1)(N+2)/6$ terms. Generating the binomials
and accumulating the output require $O(N^2)$ additional operations.
The generic combination thus costs $\Theta(N^3)$ plus evaluation of the
outer derivative vector. For sine and exponential that vector costs
$O(N)$. Peak workspace is $\Theta(N^2+M+1)$ including fixed-capacity
dual objects. These bounds count loops in fixed precision.

\subsection{Dualization of elementary functions}

As an illustrative example, the dual extension of $\sin(z)$ to order $n$
is
\begin{align}
\sin(z+\epsilon_1) &= \sum_{k=0}^n D^k(\sin z)\,\epsilon_k
= \sum_{k=0}^n \sin\!\left(z + k\frac{\pi}{2}\right)\epsilon_k.
\end{align}
Here $z$ is the scalar expansion point. To
construct the dual extension $\sin(g)$, where $g$ is a dual number, the
chain rule must be applied:
\begin{align}
\sin(g) = \sum_{k=0}^n D^k(\sin\circ g)(x_0)\,\epsilon_k.
\end{align}

Following this procedure, dual extensions for all elementary functions
can be constructed. In cases where closed-form expressions for higher
derivatives are cumbersome, functions can be dualized by combining
previously dualized elementary operations. For example, the
higher derivatives of $\arcsin(z)$ satisfy the recursive relations
\begin{align}
D^0(\arcsin z) &= \arcsin z, \\
D^1(\arcsin z) &= \frac{1}{\sqrt{1-z^2}}, \\
D^n(\arcsin z) &=
-\sum_{k=1}^{n-1} \binom{n-1}{k}
D^k\!\left(\sqrt{1-z^2}\right)
\frac{D^{n-k}(\arcsin z)}{\sqrt{1-z^2}},
\end{align}
which the FPM implementation evaluates for $n\ge 2$. It seeds a dual
variable at $z$, evaluates $\sqrt{1-z^2}$ in dual arithmetic, and
combines the resulting derivatives with lower derivatives of $\arcsin$.
The outer composition with a general dual argument again uses
\texttt{compose\_dual}. These formulas apply on a locally analytic branch.

As an example, Algorithms~\ref{LBRk} and~\ref{timesdF} show the
dualization of multiplication using the Leibniz rule.

\begin{algorithm}
\caption{Leibniz rule for the $k$-th derivative of a product.}
\label{LBRk}
\noindent\textbf{Require:} $A,B$: dualzn, $k \in \mathbb{Z}$\par
\noindent\textbf{Ensure:} $D^k(AB)$\par
\begin{enumerate}\setlength{\itemsep}{0pt}\setlength{\parskip}{0pt}
\item $res \leftarrow 0$
\item \textbf{for} $i \gets 0$ to $k$ \textbf{do}
\item \hspace*{1em}$res \leftarrow res
+ \binom{k}{i} A\%f(i) B\%f(k-i)$
\item \textbf{end for}
\item \textbf{return} $res$
\end{enumerate}
\end{algorithm}

\begin{algorithm}
\caption{Product of two dual numbers.}
\label{timesdF}
\noindent\textbf{Require:} $A,B$: dualzn\par
\noindent\textbf{Ensure:} $AB$: dualzn\par
\begin{enumerate}\setlength{\itemsep}{0pt}\setlength{\parskip}{0pt}
\item Require equal active orders $N=A\%ord=B\%ord$
\item Initialize $res$ at order $N$
\item \textbf{for} $k \gets 0$ to $N$ \textbf{do}
\item \hspace*{1em}$res\%f(k) \leftarrow \texttt{timesdzn}(A,B,k)$
\item \textbf{end for}
\item \textbf{return} $res$
\end{enumerate}
\end{algorithm}

\section{Programming with DNAOAD-FPM}
\label{sec:package}
The FPM repository provides the library, its module reference, and
runnable Fortran examples \cite{Penunuri_DNAOAD_FPM}. A standalone
distribution with an earlier interface is also available
\cite{Frpa2024}. The examples below use the FPM interface.

\subsection{Configuration and use}
The preprocessed module \texttt{config\_mod} sets the default capacity
to $M=5$ and the default precision to \texttt{real64}.
\verb|-DMAX_ORDER_DUALZN=24| selects a capacity of 24, and
\verb|-DUSE_REAL128| selects \texttt{real128} on a compiler supporting
that kind. The module \texttt{dualzn\_mod} exports the type, its
initializers, arithmetic operators, and elementary functions. The
composition routine \texttt{compose\_dual} is private.
User procedures can combine the exported operations using ordinary
Fortran expressions. The module also provides \texttt{absx}, defined
through $\sqrt{z^2}$ for complex-step applications
\cite{Squire1998,Martins2003}.

From the repository root, the following commands run two introductory
examples with GNU Fortran:
\begin{verbatim}
FPM_FC=gfortran fpm run ex1
FPM_FC=gfortran fpm run --flag "-DMAX_ORDER_DUALZN=10" ex2
\end{verbatim}
The second example evaluates $f(x)=\sin(x)\exp(-x^2)$ composed with
itself 1,000 times at $x=1.1$, through order ten. The first example
uses the default capacity and a complex expansion point
(Listing~\ref{lst:e1}). The constructor \verb|xto_dzn(z0,n)|
sets the value and order and clears active derivative components.
Seeding component one with unity makes output component $k$ equal to
the $k$th derivative. Alternatively, \verb|initialize_dualzn(r,n)|
sets the active components to zero before seeding.
\begin{lstlisting}[style=fortran-style,
caption={Derivatives of a complex power through order five.},
label={lst:e1}]
program example
  use config_mod, only: prec
  use dualzn_mod
  implicit none
  integer, parameter :: n = 5
  integer :: k
  complex(prec) :: z0
  type(dualzn) :: r, value

  z0 = (1.1_prec, 2.2_prec)
  r = xto_dzn(z0, n)
  r%f(1) = 1.0_prec
  value = sin(r)**log(r*r)
  do k = 0, n
     print *, k, value%f(k)
  end do
end program example
\end{lstlisting}

For the repeated-composition example the essential loop is
\begin{lstlisting}[style=fortran-style,caption={Core of the FPM example \texttt{ex2}, with active order ten.}]
r = xto_dzn(1.1_prec, 10)
r%f(1) = 1.0_prec
value = r
do j = 1, 1000
   value = sin(value)*exp(-value*value)
end do
\end{lstlisting}
Appendix~\ref{a1} gives shorter nested-differentiation examples in two
other languages as illustrations of their interfaces.

\subsection{Calculating gradients, Jacobians and Hessians}
The module \texttt{diff\_mod} uses univariate seeds to differentiate
functions of several variables \cite{Penunuri2024}. For a vector
function $F$, seeding $\mathbf{q}+t\mathbf{v}$ and extracting the first derivative
at $t=0$ gives $\mathbf{J}(\mathbf{q})\mathbf{v}$ in one direction. The same product can be
formed explicitly from the Jacobian. The following calls illustrate
both routes for a function of four variables with three outputs.
The repository example \texttt{ex3} also computes the gradient and
Hessian of the scalar function
$f(x,y,z)=\sin(xyz)+\cos(xyz)$.

\begin{lstlisting}[style=fortran-style,
caption={Core calls in a Jacobian--vector product calculation.
The callback \texttt{fvectest} accepts a vector of dual numbers and
returns a vector of three dual numbers.},label={lst:e3f}]
vec = [1.0_prec, 2.0_prec, 3.0_prec, 4.0_prec]
q = vec/10.0_prec + (0.0_prec, 1.0_prec)
Jmat = Jacobian(fvectest, q, 3)
JV_from_matrix = matmul(Jmat, vec)
JV_direct = d1fvector(fvectest, vec, q, 3)
\end{lstlisting}

For scalar functions, \texttt{d1fscalar} returns a first directional
derivative and \texttt{d2fscalar} returns $\mathbf{v}^{\mathsf T}\mathbf{H}\mathbf{v}$ for one
direction or $\mathbf{u}^{\mathsf T}\mathbf{H}\mathbf{v}$ by polarization for two directions.
Here $\mathbf{H}$ denotes the Hessian. For complex analytic functions these
contractions use the transpose.

\section{Cost and numerical evaluation}
\label{sec:benchmark}
We use $h(x)=\exp(\sin x)$ at $x_0=0.5$ to study the cost and
accuracy of computing all derivatives through order $N$. For this
composition, DNAOAD-FPM~\cite{Penunuri_DNAOAD_FPM} constructs two
tables of partial Bell polynomials, one for $\sin$ and one for $\exp$.
Each table supplies all derivative orders for the corresponding
function.

To measure the effect of the composition algorithm, we compared it
with an earlier implementation in DNAOAD-FPM that computed each
partial Bell polynomial $B_{n,k}$ by building a separate table.  For
each elementary function, that implementation built one table for
every pair $1\leq k\leq n\leq N$, giving $N(N+1)$ tables for the full
expression. It also computed binomial coefficients in an inner
multiplicative loop. For a requested $B_{n,k}$, the loops visited $n
k$ cells.  For rows $r\geq 3n/4$ and columns up to $\min(k,\lfloor
n/4\rfloor)$, a cell required order $n$ terms, many with binomial
coefficients that each took order $n$ multiplications. The work per
table was $\Theta(n^3 k)$, giving a total of
$\sum_{n=1}^{N}\sum_{k=1}^{n}\Theta(n^3 k)=\Theta(N^6)$.

At $N=24$, computing the full derivative vector on the same machine
took 4.57\,ms with the earlier implementation and
7.17\,$\mu$s with the shared-table implementation using GNU Fortran.
The corresponding times with Intel Fortran were 6.00\,ms and
5.70\,$\mu$s. The measured speedups within DNAOAD-FPM were 638 and
1053, respectively. For each compiler and build profile, the two
implementations agreed on the exported derivatives for both
$\exp(\sin x)$ and an exponential with a complex polynomial argument.

\subsection{Taylor control and accuracy}
In Taylor mode an elementary operation can often be propagated using a
differential identity. For example, if $y(t)=\exp(g(t))$ and $g_j,y_j$
denote normalized Taylor coefficients, matching powers in $y'=g'y$
yields $y_0=\exp(g_0)$ and
\[
 y_n=\frac{1}{n}\sum_{j=1}^n j g_j y_{n-j},\qquad n\ge1 .
\]
This rule handles any input series $g$ in quadratic time and linear
storage through order $N$. Our timing control specializes the
differential identity further for the particular composition tested.

Write $a_j=h^{(j)}(x_0)/j!$ and $c_j=\cos^{(j)}(x_0)/j!$.
Matching powers in the differential identity $h'=h\cos x$ gives
\[
 a_0=\exp(\sin x_0),\qquad
 a_{n+1}=\frac{1}{n+1}\sum_{j=0}^{n}c_j a_{n-j}.
\]
Through order $N$, this specialized Taylor recurrence requires
$\Theta(N^2)$ work and $\Theta(N)$ storage. We convert each
coefficient back to $k!a_k$ before comparing outputs. An independent
Python implementation evaluated the recurrence in decimal arithmetic
at 120 digits. Repeating at 150 digits gave identical reference
derivatives when rounded to 100 digits.

In the specialized Fortran control, cosine coefficients and factorials
were prepared once per order at the fixed point. The timed evaluation
included conversion to ordinary derivatives. At $N=24$ this control
required about 0.39\,$\mu$s with either compiler.
Both Fortran implementations used \texttt{complex(real64)}. Across the
95 pairs $(N,k)$ of requested order and derivative order, the largest
relative error for the specialized Taylor recurrence was
$1.07\times10^{-13}$ with either compiler. It occurred at $k=20$.

\subsection{Comparison with ADOL-C}
With ADOL-C 2.7.2, the C++ driver recorded $h$ on a tape and called
\texttt{hos\_forward} with a unit linear seed, real \texttt{double}
arithmetic and the tape reused. Timings include the sweep and
conversion to ordinary derivatives. A separate check on $x^3$ at $x=2$
recovered $(8,12,12,6,0)$ through order four. The first tape recording
in a separate process took 27.5\,$\mu$s (one observation). Re-recording
with the same tag took a median of 0.702\,$\mu$s. DNAOAD uses complex
double-precision coefficients with a real-valued seed, whereas ADOL-C
uses real double-precision coefficients. The timings therefore include
the cost of each scalar representation. Table~\ref{tab:audit} excludes
both tape-recording measurements.

\begin{table}[htbp]\centering\small
\caption{Median ns per complete derivative vector for $\exp(\sin x)$ at
  $x=0.5$. ADOL-C reuses a recorded tape. The two DNAOAD-FPM columns use
  the same library source compiled with GNU and Intel Fortran.}
\label{tab:audit}
\begin{tabular}{rrrr}\toprule
$N$ & ADOL-C & DNAOAD (GNU) & DNAOAD (Intel)\\\midrule
 1 &  669.1 &  144.0 &   83.6 \\
 2 &  677.4 &  173.0 &  114.5 \\
 4 &  690.3 &  257.8 &  196.4 \\
 8 &  742.1 &  607.6 &  512.8 \\
12 &  816.0 & 1335.9 & 1090.0 \\
16 &  974.4 & 2524.6 & 2048.2 \\
20 & 1096.4 & 4357.1 & 3517.5 \\
24 & 1263.3 & 7023.9 & 5567.4 \\
\bottomrule\end{tabular}\end{table}
The tests used an Intel Core i7-12650H, one logical CPU, GNU Fortran
13.3.0 (\texttt{-O3}) and Intel ifx 2025.1.1
(\texttt{-O3 -fp-model precise}). The ADOL-C driver used GNU C++ with
\texttt{-O3} and a system-installed library. After warm-up, we timed
seven batches of at least 0.12\,s each and report the median of their
mean times per evaluation. The benchmark series ran sequentially on the
same logical CPU. At the start and end of each series, the system
reported AC power and the balanced power profile. The maximum relative
errors across 95 order/component pairs were $1.11\times10^{-12}$ for
DNAOAD and $1.47\times10^{-13}$ for ADOL-C, both at derivative order 20.
DNAOAD had lower median times at the sampled orders through eight.
ADOL-C had lower median times at the sampled orders from 12 onward.
These measurements apply to this expression, point and choice of scalar
types.
\section{Programming an implicit orbital calculation in dual arithmetic}

For an elliptic orbit with semimajor axis $a$ and fixed eccentricity
$0\le e<1$, Kepler's equation $E-e\sin E=M$ determines the eccentric
anomaly $E$ from the mean anomaly $M$. The coordinates in the orbital
plane are $X=a(\cos E-e)$ and $Y=a\sqrt{1-e^2}\sin E$
\cite{JPLKepler}. Throughout this section, $M$ denotes the mean
anomaly. We use radians and set $a=1$, so that positions are
expressed in units of the semimajor axis.

We seek a local approximation to the orbital position as a function of
$M$ near a prescribed value $M_0$. We seed $\widehat M=M_0+\epsilon_1$
and initialize the constant component of $\widehat E$ with the scalar
Newton--Raphson solution at $M_0$, setting all higher components to
zero. We then perform the Newton--Raphson iteration in dual
arithmetic. For the residual $F(E,M)=E-e\sin E-M$, the derivative with
respect to $E$ gives the denominator $1-e\cos\widehat E$ in the dual
update. Seeding $\widehat M$ with unit first component makes the
components of the converged $\widehat E$ the derivatives $d^kE/dM^k$
at $M_0$. Once the iteration converges, we pass the dual eccentric
anomaly directly to both coordinate expressions and extract the
derivatives of $X$ and $Y$ with respect to $M$.

\begin{lstlisting}[style=fortran-style,caption={Newton--Raphson iteration
in dual arithmetic and evaluation of the orbital coordinates. The dual
variables \texttt{m} and \texttt{u} are initialized as described above.}]
do iter = 0, max_iter
   s = e*sin(u)
   r = u - s - m
   scale = max(1.0_prec, abs(u%f(0:n)), &
               abs(s%f(0:n)), abs(m%f(0:n)))
   ratio = abs(r%f(0:n))/scale
   if (.not. all(ieee_is_finite(ratio))) error stop
   if (maxval(ratio) <= tol) exit
   if (iter == max_iter) error stop
   delta = r / (1.0_prec - e*cos(u))
   u = u - delta
end do
x = cos(u) - e
y = sqrt(1.0_prec-e**2)*sin(u)
\end{lstlisting}

The convergence test requires every scaled residual component to be
at most $64$ times machine epsilon in double precision, with a limit
of 50 updates. Testing only the order-zero residual would accept the
initial scalar root immediately. At initialization, however, the
order-one residual is exactly $-1$ because $E_1=0$ and $M_1=1$.
Stopping at this point would give zero derivatives for both
coordinates. With the componentwise convergence test, all four
elliptic test cases required four dual updates at order 12.

With both compilers, we verified that $E=M$ for a circular orbit
($e=0$). For $e>0$, we checked the first two derivatives of $E$ with
respect to $M$ against the analytical expressions
$E'=1/(1-e\cos E)$ and $E''=-e\sin E/(1-e\cos E)^3$.
At each expansion center, we also compared the derivatives of $E$,
$X$, and $Y$ through the requested order with reference values
obtained from an independent coefficient recurrence evaluated at
110-digit precision. For each computed value $v$ and corresponding
reference value $v_{\mathrm{ref}}$, we used the scaled error
$|v-v_{\mathrm{ref}}|/\max(1,|v_{\mathrm{ref}}|)$. The largest
scaled error over all derivatives and both compilers was below
$1.52\times10^{-13}$.

\subsection{Predicting positions near the expansion center}

For each coordinate, we divided the stored derivative of order $k$
by $k!$ to obtain its Taylor coefficient and evaluated the resulting
polynomial at $\Delta M=M-M_0$ using Horner's rule. To obtain an
independent reference position at each query, we solved the scalar
equation by 325 bisection steps on $[M-e,M+e]$ in 110-digit decimal
arithmetic. We measured position error as the Euclidean distance
between the predicted and reference positions, expressed in units
of $a$.

\begin{figure}[htbp]
\centering
\includegraphics[width=9.71cm]{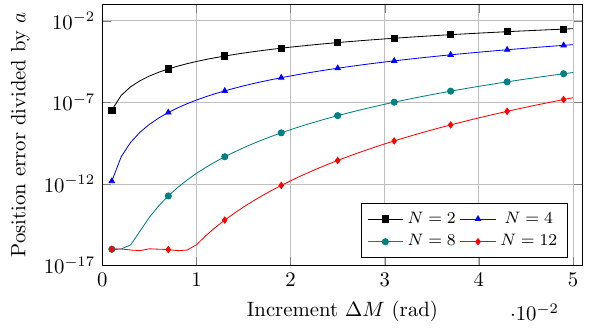}
\caption{Orbital position error for $e=0.9$ and $M_0=0.1$ using
GNU Fortran. Each curve gives the Euclidean error against the
independent reference at 50 positive increments $\Delta M$.}
\label{fig:orbit}
\end{figure}

We tested all four combinations of $e\in\{0.3,0.9\}$ and
$M_0\in\{0.1,1.0\}$ at orders $1,2,4,8,12$.
Figure~\ref{fig:orbit} shows the local position error for
$(e,M_0)=(0.9,0.1)$. At $\Delta M=0.01$, increasing the order from
8 to 12 reduced the error from $4.49\times10^{-12}a$ to
$1.89\times10^{-16}a$. At $\Delta M=0.05$, the order-12 error was
$1.93\times10^{-7}a$. The useful range of the expansion therefore
depends on the expansion center and the prescribed tolerance.

\subsection{End-to-end cost at a prescribed accuracy}

For the timing study, we selected two of the four configurations,
A with $(e,M_0)=(0.3,0.1)$ and B with $(e,M_0)=(0.9,1.0)$.
At each center, we used 256 query points $M=M_0+\Delta M$, with
equally spaced increments $\Delta M$ from 0.001 to 0.05 radians.
We applied a fixed permutation to the points and used the first $Q$
points for each group of $Q$ queries.

For each group, the dual strategy built an eighth-order expansion
once and evaluated it at all $Q$ points. The measured time included
the initial scalar solve, the dual Newton--Raphson iteration, the
calculation of both orbital coordinates, the conversion of their
derivatives to Taylor coefficients, and polynomial evaluation.
The scalar strategy solved Kepler's equation independently at each
query, using that query's mean anomaly as the initial guess for
Newton--Raphson, and then computed both coordinates.

We computed the reference positions from the exact binary values
exported by Fortran, including the rounded value of $M$ at each
query. Both strategies met the position tolerance of $10^{-12}a$ at
all 256 points at each center with both compilers. The maximum
position errors for the expansion were $3.738\times10^{-13}a$ at A
and $3.623\times10^{-14}a$ at B. The corresponding maxima for the
direct scalar solves were $1.152\times10^{-16}a$ and
$3.049\times10^{-16}a$.

\begin{table}[htbp]
\centering
\small
\caption{Median time per group of $Q=256$ queries, in
$\mu\mathrm{s}$. The ratio is the median time for independent
scalar solves divided by the median time to build and evaluate
the expansion.}
\label{tab:groups}
\begin{tabular}{llrrr}
\toprule
Center & Compiler & Scalar & Dual & Ratio \\
\midrule
A & GNU   & 16.138 & 15.660 & 1.031 \\
A & Intel & 18.605 &  9.807 & 1.897 \\
B & GNU   & 33.942 & 15.685 & 2.164 \\
B & Intel & 32.825 &  9.855 & 3.331 \\
\bottomrule
\end{tabular}
\end{table}

For each compiler, center, and group size, we measured seven paired
batches, alternating which strategy ran first. Point preparation and
reference calculations were outside the timed region. In the run
reported in Table~\ref{tab:groups}, all seven pairs favored the
expansion at $Q=256$ for each center and compiler. At A with GNU
Fortran, the expansion reduced the median group time by about
3\%. Across three runs in this configuration at $Q=256$, the
scalar-to-dual ratio ranged from 1.012 to 1.047. In the reported run,
all seven pairs favored the expansion at $Q=249$, while the paired
differences had both signs at $Q=250$ and $Q=251$.  The reported
ratios apply to scalar solves initialized independently at each
query's mean anomaly. Using the preceding scalar root as the initial
guess would constitute a different comparison.


\section{Angular kinematics of a spherical four-bar coupler}
\label{sec:coupler}
Consider the spherical four-bar linkage studied by Pe\'on-Escalante et
al.~\cite{Peon2025}. Its configuration is described by a single input
angle $\theta$. The Fortran geometry routine \texttt{fCouplerBase}
returns the three orthonormal vectors $\mathbf e_i$ of a right-handed
frame attached to the coupler as rows, with components in the fixed
frame. We use the design of Table~3 in Ref.~\cite{Peon2025},
\[
(\eta_1,\phi_1,\psi,\alpha_1,\alpha_2,\alpha_3,\alpha_4,\beta,\gamma)
=(\pi/2,0,0,\pi/2,\pi/6,\pi/2,\pi/2,0,0).
\]
The crank starts at $\theta_0=3\pi/4$ and turns at a constant
$1\,\mathrm{rad\,s^{-1}}$. With time in seconds and angles in radians,
$\theta(t)=\theta_0+t$, so a unit first component of the input angle
seeds differentiation with respect to time.

The angular velocity in fixed coordinates follows from
$\dot{\mathbf e}_i=\boldsymbol\omega\times\mathbf e_i$. With the
cyclic triples $(i,j,k)=(1,2,3),(2,3,1),(3,1,2)$, Eq.~(15) of
Ref.~\cite{Peon2025} takes the form
\begin{equation}
 \boldsymbol\omega=\sum_{\rm cyc}
 (\dot{\mathbf e}_j\cdot\mathbf e_k)\,\mathbf e_i,\qquad
 \mathbf K_n=\frac{d^{n-1}\boldsymbol\omega}{dt^{n-1}} .
 \label{eq:coupler-omega}
\end{equation}
One evaluation of the geometry at order eight supplies both operands
of each dot product. Since the dual components store ordinary time
derivatives, shifting the index by one represents $\dot{\mathbf e}_i$
through order seven. Evaluating Eq.~\eqref{eq:coupler-omega} on these
duals gives the angular velocity, acceleration, jerk and snap in
components zero through three of $\boldsymbol\omega$, and $\mathbf K_5$
through $\mathbf K_8$ in components four through seven.

\begin{lstlisting}[style=fortran-style,caption={Angular quantities from a single call to the coupler geometry routine. The second index of \texttt{e} selects the basis vector.}]
call initialize_dualzn(theta, 8)
theta(1)%f(0) = 3*pi/4
theta(1)%f(1) = 1.0_prec
base = fCouplerBase(theta, design)
call initialize_dualzn(e, 7)
call initialize_dualzn(edot, 7)
do i = 1, 3
  do j = 1, 3
    e(j,i)%f(0:7)    = base(i,j)%f(0:7)
    edot(j,i)%f(0:7) = base(i,j)%f(1:8)
  end do
end do
omega = sum(edot(:,2)*e(:,3))*e(:,1) + &
        sum(edot(:,3)*e(:,1))*e(:,2) + &
        sum(edot(:,1)*e(:,2))*e(:,3)
\end{lstlisting}

We checked the assembly of $\mathbf K_n$ from the same basis derivatives
by evaluating the Leibniz expansion in real arithmetic,
\begin{equation}
 \mathbf K_n=\sum_{\rm cyc}\ \sum_{a+b+c=n-1}
 \frac{(n-1)!}{a!\,b!\,c!}
 \bigl(\mathbf e_j^{(a+1)}\cdot\mathbf e_k^{(b)}\bigr)
 \mathbf e_i^{(c)}, \qquad a,b,c\geq 0 .
 \label{eq:coupler-leibniz}
\end{equation}
Across all eight orders, the largest absolute difference in any
component between the dual evaluation and the explicit Leibniz
expansion was $5.6\times10^{-14}$, at $n=8$.

Equation~(26) of Ref.~\cite{Peon2025} gives the angular velocity as
\begin{equation}
  \boldsymbol\omega
  = \frac{1}{2}\sum_{i=1}^{3}
    \mathbf e_i \times \dot{\mathbf e}_i .
  \label{eq:coupler-omega-cross}
\end{equation}
With the basis vectors as columns of
$\mathbf{R}=[\mathbf e_1\ \mathbf e_2\ \mathbf e_3]$, the same identity takes
the matrix form $\dot{\mathbf{R}}\mathbf{R}^{\mathsf T}=[\boldsymbol\omega]_\times$,
where $[\boldsymbol\omega]_\times\mathbf v
=\boldsymbol\omega\times\mathbf v$.
To check both the geometry and its derivatives independently of
\texttt{fCouplerBase}, we reconstructed $\mathbf{R}(t)$ from the linkage
constraints using a separate high-precision Taylor calculation.
The input crank gives
$\mathbf r_2=(\sqrt{3},\cos\theta,\sin\theta)/2$. The other moving joint
direction lies in the $xz$ plane and is perpendicular to $\mathbf r_2$.
On the chosen branch it is the unit vector along
$(-r_{2z},0,r_{2x})$. These two unit vectors give $\mathbf e_3$
and $\mathbf e_1$, respectively, with $\mathbf e_2=\mathbf e_3\times
\mathbf e_1$. From the Taylor coefficients of $\mathbf{R}$, we formed
$\mathbf{S}=\dot{\mathbf{R}}\mathbf{R}^{\mathsf T}$ by polynomial multiplication and extracted
$\boldsymbol\omega=(S_{32},S_{13},S_{21})$.
Multiplying the coefficient of degree $n-1$ by $(n-1)!$ gives the
reference vector $\mathbf{K}_n$. We computed the reference values
independently in Python using 120 decimal digits of working precision.
Repeating the calculation at 80 digits gave the same reference
components when rounded to double precision. We used the 120-digit
values to assess the double-precision Fortran results in
Table~\ref{tab:coupler}.

\begin{table}[htbp]\centering\small
\caption{Coupler angular quantities in the fixed frame at $t=0$.
Components have units of $\mathrm{rad}\,\mathrm{s}^{-n}$. The last column gives the GNU Fortran result for
$\|\mathbf K_n-\mathbf K_n^{\rm ref}\|_2/\|\mathbf K_n^{\rm ref}\|_2$
before rounding the displayed components.}
\label{tab:coupler}
\begin{tabular}{rrrrr}\toprule
$n$ & $K_{n,x}$ & $K_{n,y}$ & $K_{n,z}$ & Relative error\\\midrule
1 & 0.14285714 & 0.34992711 & -0.34992711 & $1.72\times10^{-16}$\\
2 & -0.24489796 & 0.44990628 & 0.24994793 & $2.61\times10^{-16}$\\
3 & -0.13994169 & -0.09283780 & 0.49275450 & $5.71\times10^{-16}$\\
4 & 0.85964182 & -0.82941906 & 0.47235058 & $3.56\times10^{-16}$\\
5 & 2.10198132 & -2.65469481 & -0.95271699 & $5.73\times10^{-16}$\\
6 & 0.68379672 & -0.78182381 & -8.24224389 & $1.14\times10^{-15}$\\
7 & -25.19393401 & 29.60823995 & -17.72050389 & $9.90\times10^{-16}$\\
8 & -132.22805367 & 149.74324894 & 56.02084344 & $7.36\times10^{-16}$\\
\bottomrule\end{tabular}
\end{table}

All twelve components through snap reproduce Table~4 of
Ref.~\cite{Peon2025} when rounded to its four decimal places.
For the four higher orders, Table~\ref{tab:coupler} gives relative
errors below $1.2\times10^{-15}$ against the independent reference.
With both GNU Fortran and Intel Fortran, the \texttt{spherical\_4r}
example reproduces all 24 Cartesian components in
Table~\ref{tab:coupler} to the eight decimal places shown.

Starting from the seeded input angle, the complete calculation of the
24 Cartesian components in Table~\ref{tab:coupler} took
$105.41\,\mu\mathrm{s}$ per call with GNU Fortran and
$99.23\,\mu\mathrm{s}$ with Intel Fortran on the Intel Core i7-12650H.
Each call evaluated the coupler frame through order eight and returned
$\mathbf K_1,\ldots,\mathbf K_8$ in fixed coordinates. Evaluating the
frame alone took $104.28$ and $96.63\,\mu\mathrm{s}$, respectively.
After warm-up, we timed seven batches per evaluator, alternating their
order, and report the median of the seven mean times per call.

\section{Conclusions}\label{conclusions}
DNAOAD-FPM allows a Fortran routine to return a numerical result
together with its ordinary derivatives through a selected order.
The calling routine can use that result directly in further
calculations. The two applications exploit this property in different
ways: the orbital solver passes an implicitly defined quantity to a
coordinate calculation, and the linkage geometry supplies a moving
frame from which angular kinematic quantities are obtained.

The cost of this approach depends in part on how elementary functions
are evaluated. The generic composition routine builds one table of
partial Bell polynomials for all requested orders. Given the
derivatives of the outer function, combining them with those of the
inner function requires cubic work in the active order. A Taylor
recurrence specialized to $\exp(\sin x)$ requires quadratic work.
For this expression at $x=0.5$, DNAOAD-FPM had lower median
evaluation times at the tested orders through eight, while ADOL-C's
higher-order forward sweep had lower medians at the tested orders
from 12 to 24. These timings include complex double-precision
arithmetic in DNAOAD-FPM and real double-precision arithmetic with
a reused tape in ADOL-C.

Kepler's equation illustrates a consequence of using a dual result
in a subsequent calculation. When the Newton--Raphson iteration
started at the converged scalar root, the order-zero residual already
met the tolerance, although its first derivative with respect to
mean anomaly was still $-1$. Testing every active residual component
allowed the iteration to recover the derivatives needed for the
orbital coordinates. At each of the two timing centers, a single
eighth-order expansion then met the position tolerance of $10^{-12}a$
at all 256 queries. Including construction, the expansion took about
as long as independent scalar solves at center A with GNU Fortran.
In the other three center--compiler combinations, the scalar solves
took 1.90 to 3.33 times as long. The usefulness of the expansion
therefore depends on how many queries share a center and on the
interval over which it meets the required accuracy.

For the spherical four-bar linkage, one order-eight evaluation of
the existing geometry routine supplied the coupler frame and its
time derivatives. Since the dual components store ordinary
derivatives, shifting them by one position gave the time derivative
of each frame vector through order seven. The same angular velocity
expression then yielded the coupler's angular velocity and its first
seven time derivatives in fixed coordinates. At the tested
configuration, the results through angular snap reproduced the
published components to four decimal places, and the four higher
orders agreed with an independent reconstruction from the linkage
constraints to within $1.2\times10^{-15}$ in relative Euclidean norm.
Thus, extending the calculation beyond snap required derivatives of
the frame through the corresponding order, while the geometry routine
and the expression for angular velocity served all eight orders.

\section*{Code availability}
The DNAOAD-FPM library and the executable spherical four-bar example
(\texttt{spherical\_4r}) are available in the project
repository~\cite{Penunuri_DNAOAD_FPM}.

\appendix
\section{Examples of repeated differentiation}\label{a1}
The following listings compute derivatives through order ten for five
compositions of $f(x)=\sin(x)\exp(-x^2)$ at $x=1.1$. They illustrate
recursive differentiation and graph construction. Table~\ref{tab:audit}
compares DNAOAD and ADOL-C on a different test function.

\begin{lstlisting}[style=python-style,caption={Repeated differentiation in PyTorch.}]
import torch

def ftest(x):
    for _ in range(5):
        x = torch.sin(x) * torch.exp(-x*x)
    return x

x = torch.tensor(1.1, dtype=torch.float64, requires_grad=True)
y = ftest(x)
print(0, y.item())
for k in range(1, 11):
    y, = torch.autograd.grad(y, x, create_graph=(k < 10))
    print(k, y.item())
\end{lstlisting}

\begin{lstlisting}[style=julia-style,caption={Nested differentiation in ForwardDiff.}]
using ForwardDiff

function ftest(x)
    for _ in 1:5
        x = sin(x) * exp(-x*x)
    end
    return x
end

function nth_derivative(f, x, n)
    n == 0 && return f(x)
    return ForwardDiff.derivative(y -> nth_derivative(f, y, n-1), x)
end

for k in 0:10
    println(k, " ", nth_derivative(ftest, 1.1, k))
end
\end{lstlisting}


\end{document}